\documentclass[12pt,a4paper]{article}
\usepackage{amssymb}
\usepackage{amsthm}
\usepackage[utf8]{inputenc}
\usepackage[english]{babel}
\usepackage{color}
\usepackage{amsmath}
\usepackage{enumerate}

\newtheorem{Thm}{Theorem}[section]
\newtheorem{Lem}[Thm]{Lemma}
\newtheorem{Prop}[Thm]{Proposition}
\newtheorem{Cor}[Thm]{Corollary}

\theoremstyle{Def}

\theoremstyle{Rem}

\newcommand{\C}{\mathbb{C}}
\newcommand{\Hi}{\mathcal{H}}
\newcommand{\R}{\mathbb{R}}

\newcommand{\tk}{\mathcal{T}_{q+1}}
\newcommand{\bordtk}{\partial\mathcal{T}_{q+1}}

\title{Functions conditionally of negative type on groups acting on regular trees}
\author{Antoine GOURNAY\footnote{Supported by the ERC-StG 277728 ``GeomAnGroup''.}~~and Pierre-Nicolas JOLISSAINT\footnote{Supported by Swiss SNF project 20-137696.}}

\begin{document}

\maketitle

\begin{abstract} 
Let $\tk=(V,E)$ be the $(q+1)$-regular tree and let $G$ be a group of automorphisms acting transitively on the vertices and on the boundary of $\tk$. We give an upper bound for the growth of cocycles with values in any unitary representation of the group $G$. This bound is optimal by projecting the Haagerup cocycle onto an appropriate subspace of $\ell^{2}(E)$. We also obtain a description of functions conditionally of negative type which are unbounded.\\
\end{abstract}


\section{Introduction}

A second countable, locally compact group $G$ has the Haagerup property if it admits a proper affine isometric action on some Hilbert space $\Hi$. It can be seen both as a generalisation of amenability and as a strong negation of Property (T). Its importance is due, at least in part, to its connection with the Baum-Connes Conjecture, since it holds for groups having the Haagerup property (see \cite{KY} and \cite{Tu}). The class of groups having this property is large and the most famous examples are : amenable groups, free groups and, more generally, groups acting properly on trees. In order to have a better understanding of these groups as well as to single out interesting representations, a lot of efforts has been dedicated to investigate their unitary dual. The starting point of our work was the following theorem of C. Nebbia about the $1$-cohomology of some large groups acting on regular trees. In the rest of this text, for an integer $q\geq 2$, we will denote by $\tk$ the homogeneous $(q+1)$-regular tree and by $\bordtk$ its boundary.

\begin{Thm}\label{ThmNeb}(\cite{N})
Let $G$ be a closed non-compact subgroup of $\mathrm{Aut}(\tk)$, with $q\geq 2$. Suppose that $G$ acts transitively on the vertices and on the boundary of $\tk$. Then, there exists a unique irreducible representation $\sigma^{-}$ such that $H^{1}(G,\sigma^{-})\neq 0$. Furthermore, the vector space $H^{1}(G,\sigma^{-})$ is of (complex) dimension $1$. 
\end{Thm}

This result implies that, up to some renormalisations, there exists a unique unbounded cocycle $b$ relatively to the representation $\sigma^{-}$. The present paper is devoted to identify this cocycle. One way of characterising $b$ is by the computation of its associated function conditionally of negative type, that is, $g \mapsto \|b(g)\|^{2}$. In order to estimate such quantities, we will exploit the extra structure given by the graph on which the group acts to define a notion of harmonicity for cocycles. Let $G$ be a topological group acting transitively on a graph with bounded degree $X=(V,E)$. Fix a basepoint $x_{0}\in V$ and suppose that the subgroup $K=\textrm{Stab}(x_{0})$ is compact in $G$. By transitivity, we obtain a natural identification between $V$ and $G/K$. In particular, if $b : G \rightarrow \Hi$ is a cocycle taking value in some Hilbert space $\Hi$, then it is possible to translate this cocycle in such a way that it vanishes on $K$. Such a cocycle factors through $G/K$ and the corresponding map can be seen as a $G$-equivariant embedding $F : V \rightarrow \Hi$. This transformation allows us to consider the action of the Laplace operator. Recall that for a locally finite, non-oriented graph $X=(V,E)$, we will denote by $\mathcal{L}$ the \textbf{Laplace operator} defined as
$$
\mathcal{L}f(x)=\left(\frac{1}{\mathrm{deg}(x)}\sum_{y\sim x} f(y)\right)-f(x),
$$
for any function $f : V \rightarrow \C$ and any vertex $x\in V$. More generally, this defines an operator acting on maps taking value in any Hilbert space $\Hi$. It is therefore natural to call a map $F : V \rightarrow \Hi$ \textbf{harmonic} if $\mathcal{L}F$ is identically $0$. Going back to a cocycle $b$ which is identically $0$ on $K$, we say that $b : G \rightarrow \Hi$ is harmonic if and only if it factors through a harmonic map $F : V \rightarrow \Hi$.\\ 
The main result of our paper is the following theorem which quantifies the growth of $G$-equivariant embeddings of $V$ into some Hilbert space.

\begin{Thm}\label{compression}
Let $G$ be a closed non-compact subgroup of $\mathrm{Aut}(\tk)$, with $q\geq 2$. Suppose that $G$ acts transitively on the vertices and on the boundary of $\tk$. Assume furthermore that $G$ acts on a Hilbert space $\Hi$ by affine isometries. Then : 

\begin{enumerate}\renewcommand{\labelenumi}{(\roman{enumi})}
	\item Any $G$-equivariant map $F : V \rightarrow \Hi$ such that $F(x_{0})=0$ for some $x_{0}\in V$ satisfies:
	$$
	\|F(x)\|^{2}\leq Ad(x,x_{0})-B+Bq^{-d(x,x_{0})},
	$$
	where
	$$
	A=\frac{(q+1)\|F(x_{1})\|^{2}}{q-1}, B=\frac{2q\|F(x_{1})\|^{2}}{(q-1)^{2}}
	$$
	and $x_{1}$ is any vertex adjacent to $x_{0}$.
	\item Furthermore, if the map $F$ here above is harmonic and non-constant, then the equality occurs.
	\item There exists a $G$-equivariant map $F : V \rightarrow \Hi$ which is harmonic and non-constant.
\end{enumerate}
\end{Thm}

Therefore, the upper bound appearing in the first claim here above is optimal.\\
Combining this with Nebbia's result, we deduce that this map is essentially unique and this allows us to give a complete description of functions conditionally of negative type on groups acting transitively both on $\tk$ and $\bordtk$. Interestingly, the optimal cocycle appearing in Theorem \ref{compression} can be concretely realised as the projection of the Haagerup cocycle onto an appropriate $G$-invariant closed subspace of $\ell_{\mathrm{alt}}^{2}(\mathbb{E})$, where $\mathbb{E}$ is the set of oriented edges of $\tk$.\\
Another motivation for the present paper is related to some questions about coarse (equivariant) embeddings of groups into Hilbert spaces. Roughly speaking, the problem is to represent a group $G$ together with a metric $d_{G}$ into a Hilbert space $\Hi$ by preserving large distances. More precisely, a map $F:G \to \Hi$ is called a \textbf{coarse embedding} if we can find control functions $\rho_\pm:\R^{+} \to \R^{+}$ which satisfy $\lim_{t\rightarrow +\infty }\rho_{\pm}(t)=+\infty$ and
$$
\rho_-(d_{G}(x,y)) \leq \|F(x)-F(y)\| \leq \rho_+(d_{G}(x,y)).
$$
for all $x,y\in G$. When the embedding $F$ is equivariant, it turns out that $F$ is (up to some translation) a cocycle $b$. In general, it is difficult to find examples of coarse embeddings together with a precise description of their control functions. From this point of view, Theorem \ref{compression} as well as various heuristics indicate that harmonic cocycles provide optimal embeddings (for example, Remark 2.6 in \cite{NP} shows that such cocycles are optimal in expectation for the distribution of the random walk). 
Actually, a result of Guichardet says that we can always restrict our attention to harmonic cocycles, when considering coarse equivariant embeddings of groups. Before stating this result in our setting, let us introduce two more definitions. Let $G$ be a topological group, let $K$ be a compact subgroup and let $\pi$ be a unitary representation of $G$ acting on $\Hi$. We define $Z_{K}^{1}(G,\pi)$ (resp. $B_{K}^{1}(G,\pi)$) to be the space of cocycles (resp. coboundaries) with respect to $\pi$ which are identically $0$ on $K$. Recall that such cocycles factor through $G/K$. 

\begin{Prop}\label{PropGuich}(Théorème 3, \cite{Gu})
Let $G$ be a topological group acting transitively on a graph $X=(V,E)$ with bounded degree. Suppose that the subgroup $K=\textrm{Stab}(x_{0})$ is compact in $G$. Let $\pi$ be a unitary representation of $G$ without non-trivial invariant vectors. Then, endowing $Z_{K}^{1}$ with the topology of uniform convergence on compact sets, we obtain : 
$$
Z_{K}^{1}(G,\pi)=(\ker \mathcal{L}\cap Z_{K}^{1}(G,\pi))\oplus\overline{B_{K}^{1}(G,\pi)}.
$$
In particular, the reduced cohomology group is given by :
$$
\overline{H}^{1}(G,\pi)=\ker \mathcal{L}\cap Z_{K}^{1}(G,\pi).
$$  
\end{Prop}

\smallskip

The following direct corollary implies that the behaviour at infinity of control functions of $\pi$-equivariant embeddings can be determined by analysing the growth of functions $f : G \rightarrow \R^{+}$ which are non-harmonic but whose Laplacian is constant.


\begin{Cor}\label{CorCompression}
Let $G$ be a discrete group endowed with a length function denoted by $|\cdot|$. Let $\pi$ be a unitary representation without almost invariant vectors and of non-trivial cohomology. Let $\mathcal{F}$ be the set of all functions $f:G \to \mathbb{R}$ such that $f \geq 0$, $f(e) = 0$ and $\mathcal{L} f \equiv 1$. If $\rho_-(t)$ is the lower control function of a $\pi$-equivariant embedding, then there exists a function $f \in \mathcal{F}$ so that $\rho_{-}(t)$ is bounded (above and below, up to multiplicative and additive constants) by $\min \{ \sqrt{f(g)} \ : \ |g| = t \}$, uniformly on $t$.
\end{Cor}

A second consequence of Guichardet's result is related to the reduced cohomology of groups with infinitely many finite conjugacy classes. 

\begin{Prop}\label{PasHarmRepC0}
Let $G$ be a finitely generated group that admits infinitely many finite conjugacy classes. Then, for any $C_{0}$-representation $\pi$ of $G$, we have :
$$
\overline{H^{1}}(G,\pi)=0,
$$
\end{Prop}

\medskip

The paper is organised as follows : We start with the decomposition of $Z_{K}^{1}$ in terms of harmonic cocycles in Section 2. As a consequence, we then prove Corollary \ref{CorCompression} in \S{}\ref{ConsequencesOfDecompo1} and Proposition \ref{PasHarmRepC0} in \S{}\ref{RepC0harm}. In Section 3, we specialise to groups acting on trees. Theorem \ref{compression}.(i)-(ii) is proved in \S{}\ref{secrecrel} and Theorem \ref{compression}.(iii) in \S{}\ref{secequivharm}. Section 4 is dedicated to the decomposition of the Haagerup cocycle. We single out its harmonic part by projecting it onto an appropriate subspace of $\ell_{\mathrm{alt}}^{2}(\mathbb{E})$ and, in \S{}\ref{VirtualCoboundaries}, we show that it can be described as a virtual coboundary. In Section 5, we discuss functions conditionally of negative type and provide various examples. As an application of Theorem \ref{compression}, we give a comprehensive description of pure elements of $\mathrm{CL}(G)$, the convex, positive cone of negative type functions on $G$.
\medskip

{\bf Acknowledgements:} We are grateful to Alain Valette to have suggested to look at the projection of the Haagerup cocycle and to have shown Lemma \ref{NegaLap} to us.

\section{Cohomology and Harmonicity}

\subsection{Definition of cohomology}

We first fix some notations concerning cocycles and affine isometric actions (we refer to Chapter 2 of \cite{BHV} for a complete discussion). Let $G$ be a second countable, locally compact group and let $\pi$ be a continuous unitary representation on a Hilbert space $\Hi$. A \textbf{cocycle} with values in $\pi$ is a continuous function $b : G \rightarrow \Hi$ satisfying the so-called cocycle relation, that is, $b(gh)=\pi(g)b(h)+b(g)$, for every $g,h\in G$. A cocycle of the form $g \mapsto (\pi(g)-1)v$, for some vector $v\in\Hi$, is called a \textbf{coboundary}. We denote by $Z^{1}(G,\pi)$ (resp. $B^{1}(G,\pi)$) the space of cocycles (resp. coboundaries) with values in $\pi$. The first cohomology space of $\pi$ is defined as the quotient space 
$$
H^{1}(G,\pi)=Z^{1}(G,\pi)/B(G,\pi).
$$ 
For later use, we introduce the so-called \textbf{coboundary map} $d : \Hi \rightarrow Z^{1}(G,\pi)$ defined by $(d\xi)(g)=(\pi(g)-1)\xi$, for any $g\in G$ and $\xi\in\Hi$. The space of coboundaries corresponds to the image of $\Hi$ via $d$.\\
Note that, a map $b$ satisfies the cocycle relation, if and only if, $G$ acts by affine isometries on $\Hi$ by $\alpha(g)v=\pi(g)v+b(g)$. It is clear from the definition that such an affine action $\alpha$ has a fixed point if and only if $b$ is a coboundary. We will say that $\pi$ (resp. $b$) is the linear (resp. translation) part of the affine isometric action $\alpha$, and use the self-explanatory notation $\alpha=(\pi,b)$.\\
The proofs of this section (with the exception of the  \S{}\ref{RepC0harm}) may be found in \cite[\S{}3-\S{}5]{Gu}.

\subsection{Link between $Z_{K}^{1}$ and $Z^{1}$}

As explained in the introduction, given a unitary representation $\pi$, it is sometimes more convenient to work with the space of cocycles vanishing identically on some compact subgroup $K$ of $G$, than to work with the whole space $Z^{1}(G,\pi)$. The goal of this subsection is to give a precise link between the space $Z^{1}(G,\pi)$ and its subspace $Z_{K}^{1}(G,\pi)$ formed by all the cocycles which are identically $0$ on $K$.\\    
Let $K$ be a compact subgroup of $G$. For any continuous map $f : G \rightarrow \Hi$, one can define its average on $K$ by
$$
A_{K}f=\int_{K} f(k) \, dk,
$$
where $dk$ represents the Haar measure of $G$, normalised so that $K$ has measure $1$. Taking a cocycle $b\in Z^{1}(G,\pi)$ and using the fact that $dk$ is invariant, it is easy to check that the cocycle $b+d(Ab)$ vanishes identically on $K$. It is worth noting that this procedure does not change the cohomology class of $b$ in $Z^{1}(G,\pi)$. Therefore, it is natural to introduce $Z_{K}^{1}(G,\pi)$ (resp. $B_{K}^{1}(G,\pi)$) the subspace of $Z^{1}(G,\pi)$ formed by the cocycles (resp. coboundaries) which are identically $0$ on $K$. Note that $B_{K}^{1}$ is the image via the coboundary map $d$ of the subspace of $\Hi$ formed by $\pi(K)$-fixed vectors. The following well-known lemma gives the relation between those spaces and states that we do not lose any information on the cohomology level.
\begin{Lem}\label{HKisH}
Let $\Psi : Z^{1}(G,\pi) \rightarrow Z_{K}^{1}(G,\pi)$ be the map defined by
$$
\Psi(b)=b+d\left(A_{K}b\right).
$$
Then, the map $\Psi$ induces an isomorphism between the first cohomology group $H^{1}(G,\pi)$ and the quotient space $Z_{K}^{1}(G,\pi)/B_{K}^{1}(G,\pi)$.
\end{Lem}

We stress that the average argument in the result here above already appears in Section 2 of \cite{N}. 

%
%

\subsection{Hilbert space structure on $Z_{K}^{1}(G,\pi)$}\label{SubdecompositionCoho}

In this subsection, under certain assumptions on the group $G$, we will set a Hilbert space structure on $Z_{K}^{1}(G,\pi)$ in order to decompose it in a natural way. For the rest of the section, we will assume that $G$ is a locally compact group which is generated by some subset $K\cup S$, where $K$ is a compact subgroup and $S$ is a finite set. This means, for instance, that any element $g\in G$ can be expressed as a product $g=\prod_{j=1}^{n}s_{j}k_{j}$, with $s_{j}\in S$ and $k_{j}\in K$ for all $j$.\\
The reader should have the following two examples in mind.
\begin{enumerate}
	\item Let $G$ be a finitely generated group. Take $S$ to be a finite generating set and take $K$ to be the trivial group.\\ 
	\item The second example is given by a topological group $G$ acting transitively on some regular graph $X=(V,E)$ of finite degree. We assume that the topology on $G$ is such that the stabiliser of any vertex is compact. In this case, we fix a base-vertex $x_{0}\in V$ and we let $K$ to be the vertex stabiliser of $x_{0}$. Choose $S$ so that $S\cdot x_{0}$ is the $1$-sphere centred in $x_{0}$ and $|S|$ is the degree of $X$. By Milnor-Schwarz Lemma, the group $G$ is generated by the compact set $K\sqcup S$.
\end{enumerate} 
Thanks to the cocycle relation, a cocycle $b$ is completely determined by the values $\{b(h)\}_{h\in K\cup S}$. In particular, since the elements of $Z_{K}^{1}(G,\pi)$ are exactly the cocycles which are constant on cosets of the form $gK$, it is natural to introduce the following scalar product
$$
(b,b')_{Z^{1}_{K}}=\frac{1}{|S|}\sum_{s\in S} \langle b(s),b'(s)\rangle,
$$ 
for $b,b'\in Z_{K}^{1}$. Since $b$ is identically $0$ on $G$ if and only if $b$ is identically $0$ on $S$, this scalar product is non-degenerate. We write $\|\cdot\|_{Z_{K}^{1}}$ for the induced norm on $Z_{K}^{1}$. To be sure that this scalar product defines a Hilbert space, we need to check that the space $Z_{K}^{1}(G,\pi)$ is complete for the norm $\|\cdot\|_{Z_{K}^{1}}$. 
However, the classical topology on $Z^{1}(G,\pi)$ is given by the topology of uniform convergence on compact subsets of $G$. Since under our assumptions $G$ is $\sigma$-compact, it is well-known that this topology turns $Z^{1}(G,\pi)$, and therefore $Z_{K}^{1}(G,\pi)$, into a Fr\'echet space. It is easy to show that it again gives the same topology. This is the content of the next lemma (see \cite[\S{}4]{Gu}).\\

\begin{Lem}\label{LemCompleteness}
The topology on $Z_{K}^{1}(G,\pi)$ induced by the norm $\|\cdot\|_{Z_{K}^{1}}$ is equivalent to the topology of uniform convergence on compact sets.
\end{Lem}

\subsection{Orthogonal decomposition of $Z_{K}^{1}(G,\pi)$}

To decompose $Z_{K}^{1}$, we start with the coboundary map. Let $d_{K} : \Hi^{\pi(K)} \rightarrow Z_{K}^{1}(G,\pi)$ be the restriction of the coboundary map to the subspace of $\pi(K)$-invariant vectors. Recall that its image corresponds to $B_{K}^{1}$. 
Simple computations give:
$$
d_{K}^{\ast}b=-\frac{1}{|S|}\sum_{s\in S}(b(s^{-1})+b(s)).
$$ 
Going back to the decomposition, we trivially have $Z_{K}^{1}(G,\pi)=\ker d_{K}^{\ast} \oplus (\ker d_{K}^{\ast})^{\bot}$. It is a general fact that the orthogonal complement of $\ker d_{K}^{\ast}$ coincides with $\overline{\textrm{im } d_{K}}$. Moreover, this latter space is just the closure of $B_{K}^{1}(G,\pi)$ with respect to the norm $\|\cdot\|_{Z_{K}^{1}}$. Therefore, we get the following general orthogonal decomposition.

\begin{Prop}\label{PropGenDecompo}
Let $G$ be a topological group which is generated by a set $K\cup S$, where $K$ is a compact subgroup and $S$ is a finite set. Then, we have :
$$
Z_{K}^{1}(G,\pi)=\ker d_{K}^{\ast}\oplus\overline{B_{K}^{1}(G,\pi)}.
$$
\end{Prop}

We define $\ker d_{K}^{\ast}$ to be the space of \textbf{harmonic cocycles}. This term has the following explanation.\\
Let $X=(V,E)$ be the Schreier graph defined as follows. The set of vertices $V$ is identified with $G/K$ and (oriented) edges are of the form $(gK,gsK)$, for $g\in G$ and $s\in S$. We note that the group naturally acts on $X$ by left multiplication. Up to enlarging $S$, we can suppose that for any $s\in S$ there exists a unique $s'\in S$ such that $ss'\in K$. This condition is slightly weaker than assuming $S$ to be symmetric (i.e. $S=S^{-1}$) and it ensures that the graph $X$ is non-oriented. One can express the Laplacian of a map $F : V \rightarrow \Hi$ by $\mathcal{L}F(gK)=-F(gK)+\frac{1}{|S|}\sum_{s\in S} F(gsK)$, for any $g\in G$. In particular, let $b\in Z_{K}^{1}(G,\pi)$ be a cocycle and denote by $F : V \rightarrow \Hi$ the map $b$ seen as a function on the graph $X$. Then, we claim that the map $F$ is harmonic, namely, $F\in\ker\mathcal{L}$, if and only if $b$ belongs to $\ker d_{K}^{\ast}$. Firstly, we remark that, although $S$ is not necessarily symmetric, to any coset $s^{-1}K$ with $s\in S$, there exists a unique corresponding coset $s'K$ with $s'\in S$. Therefore, we clearly have $\sum_{s\in S}b(s^{-1})=\sum_{s\in S}b(s)$. Hence, we may write :
$$
d_{K}^{\ast}b=-\frac{2}{|S|}\sum_{s\in S}b(s).
$$
Secondly, since the map $F$ obtained from $b$ is $G$-equivariant, that is $F(gK)=\alpha(g)F(K)$, for all $g\in G$, where $\alpha=(\pi,b)$, it is easy to see that
$$
\mathcal{L}F(gK)=\pi(g)\mathcal{L}F(K).
$$ 
This implies that $F$ is harmonic if and only if $F$ is harmonic at the vertex represented by $K$. Finally, we obtain the relation $\mathcal{L}F(K)=-b(e)+\frac{1}{|S|}\sum_{s\in S}b(s)=-\frac{1}{2}d_{K}^{\ast}b$, which implies the claim.\\ 
It is a good time to identify the different objects we defined in the case of the two examples cited at the beginning of Subsection \ref{SubdecompositionCoho}.
\begin{enumerate}
	\item If $G$ is a finitely generated group with finite symmetric generating set $S$ and compact subgroup $K=\{e\}$, then $\Hi^{\pi(K)}=\Hi$, $Z_{K}^{1}=Z^{1}$ and the graph $X$ defines here above is the Cayley graph of $G$ with respect to $S$. These facts imply that the decomposition of Proposition \ref{PropGenDecompo} holds for $Z^{1}$ and that the harmonicity of a cocycle can be read off directly on the Cayley graph of $G$.\\
	\item In the second example, we note that the graph we built here above coincides with the graph on which $G$ acts.
\end{enumerate}

\subsection{Some consequences of Proposition \ref{PropGenDecompo}}\label{ConsequencesOfDecompo1}

A first straightforward consequence is the description of the reduced cohomology in terms of harmonic cocycles. Recall that the \textbf{reduced cohomology} group of $G$ taking value in $\pi$ is defined as the quotient space
$$
\overline{H^{1}}(G,\pi)=Z^{1}(G,\pi)/\overline{B^{1}(G,\pi)},
$$
where $\overline{B^{1}(G,\pi)}$ is the closure of the space of coboundaries with respect to the topology of uniform convergence on compact sets.
 
\begin{Cor}\label{CorRedCoh}
We have the isomorphism :
$$
\overline{H^{1}}(G,\pi)=\ker d_{K}^{\ast}.
$$
In particular, if the representation $\pi$ does not have almost invariant vectors, then
$$
H^{1}(G,\pi)=\ker d_{K}^{\ast}.
$$
\end{Cor}
\medskip

In the case where the subspace $\Hi^{\pi(K)}$ reduces to $0$, then the space of coboundaries $B_{K}^{1}$ reduces to the constant cocycle $0$. Hence, we get :

\begin{Cor}\label{CorNonSphericalHarmonic}
If the representation $\pi$ has no non-trivial $K$-invariant vectors, then
$$
H^{1}(G,\pi)=Z_{K}^{1}(G,\pi)=\ker d_{K}^{\ast}.
$$
That is, any cocycle vanishing identically on $K$ is harmonic.
\end{Cor}

This result will be useful when we will discuss non-spherical representations of closed subgroups of $\textrm{Aut}(\tk)$. Indeed, it is straightforward that for such a representation $\pi$, the cohomology group $H^{1}(G,\pi)$ is trivial if and only if every cocycle which is identically $0$ on $K$ is trivial, as noticed by Nebbia in Section 2 of \cite{N}. \\
We are ready to prove Corollary \ref{CorCompression}.

\noindent\textbf{Proof of Corollary \ref{CorCompression} : } If a cocycle $b$ is harmonic, a simple computation (see Lemma \ref{lapandgradient}) shows that the function $f(g) = \|b(g)\|^2$ will satisfy, $f \geq 0$, $f(e)=0$ and $\mathcal{L}f \equiv c$ where $c = \|b\|_{Z^1_K}^2$. Up to a scalar multiplication of $b$, one may normalise this constant $c$ to $1$.  Next, because $\pi$ does not contain almost invariant vectors, any cocycle $b'$ may be written as $b' = cb + z$ where $z$ is a coboundary and $c \in \R^+$ accounts for the fact that we normalised $\mathcal{L}\|b(\cdot)\|^2\equiv1$. Since $\|z(g)\| = \|\pi(g) \xi - \xi\|$ is bounded (uniformly) by $2\|\xi\|$, $c\|b(g)\| - 2\|\xi\| \leq  \|b'(g)\| \leq c\|b(g)\| + 2\|\xi\|$. Replacing $\|b(g)\|$ by $\sqrt{f(g)}$ yields the claim. \hfill $\square$

A last consequence of Proposition \ref{PropGenDecompo} concerns harmonic functions on graphs whose gradient tends to $0$ at infinity. More precisely, let $G$ be a finitely generated group with finite symmetric generating set $S$ and let $X$ be the Cayley graph of $G$ with respect to $S$. If $b\in Z^{1}(G,\pi)$ is an harmonic cocycle, then one can define, for any $\xi\in\Hi$, an harmonic function $\varphi_{\xi} : G \rightarrow \C$ by $\varphi_{\xi}(g)=\langle b(g),\xi\rangle$.  Let us check it is indeed harmonic. Using the cocycle relation for $b$, we get :
\begin{eqnarray*}
\sum_{s\in S} \varphi_{\xi}(gs)
&=&
\sum_{s\in S}\langle b(gs),\xi\rangle \\
&=&
|S|\langle b(g),\xi\rangle + \langle \sum_{s\in S} b(s),\pi(g^{-1})\xi\rangle \\
&=&
|S|\varphi_{\xi}(g),
\end{eqnarray*}    
where the last equality follows from the harmonicity of $b$. Furthermore, the gradient of $\varphi_{\xi}$ is given by 
$$
\nabla\varphi_{\xi}(g,gs)=\varphi_{\xi}(gs)-\varphi_{\xi}(g)=\langle b(s),\pi(g^{-1}\xi)\rangle,
$$
for any edge $(g,gs)$. In particular, if the representation $\pi$ is $C_{0}$, then the gradient $\nabla\varphi_{\xi}$ is a $C_{0}$ function on the set of edges\footnote{Recall that a function $f : W \rightarrow \C$ on some countable set $W$ is $C_{0}$ if, for any $\epsilon >0$, we can find a finite set $F\subset W$ such that $|f|$ is strictly smaller than $\epsilon$ outside of $F$. A representation of a countable group $G$ is $C_{0}$ if all of its coefficients, that is the maps of the form $g \mapsto \langle\pi(g)x,y\rangle$, for $x,y\in\Hi$, are $C_{0}$ functions.}. In particular, we get the corollary :

\begin{Cor}
Let $G$ be a finitely generated group and let $S$ be a finite symmetric generating set. If the Cayley graph $X$ of $G$ with respect to $S$ has no non-constant harmonic function with $C_{0}$ gradient, then
$$
\overline{H^{1}}(G,\pi)=0,
$$
for any $C_{0}$ representation $\pi$.
\end{Cor}

\subsection{Harmonic functions with gradient in $C_0$}\label{RepC0harm}

It is fairly easy to check that a group which admits a harmonic function with gradient in $C_0$ has infinitely many Lipschitz harmonic functions (hence may not be nilpotent, see \cite{CM} or \cite[Theorem 1.4]{K}). An argument, which may be found in \cite[\S{}2.2]{MY}, also shows that the lamplighter (with finite lamp states) on $\mathbb{Z}$ has no sublinear harmonic functions (hence no harmonic function with gradient in $C_0$). G.~Kozma pointed out to the first named author that this is also the case for solvable Baumslag-Solitar groups. The same is conjectured to hold for polycyclic groups.

\begin{Lem}
If $G$ has infinitely many finite conjugacy classes, then there are no harmonic functions with gradient in $C_0$.
\end{Lem}
It is probably nicer to consider the particular of a group with an infinite centre when reading the upcoming proof for the first time.
\begin{proof}
Define a transport pattern (see also \cite{Go}) from $\phi$ to $\xi$ (two finitely supported measures) to be a finitely supported function on the edges $\tau$ so that $\nabla^* \tau = \xi - \phi$. Next, note that if $h$ is harmonic, then $\langle h \mid P^n_g \rangle = h(g)$ where $\langle h \mid f \rangle = \sum_{x \in X} h(x)f(x)$ (if at least one of $h$ or $f$ has finite support) and $P^n_g$ is the distribution of the simple random walk starting at $g$ at time $n$. Hence, if $C$ is a finite conjugacy class,
\[
\frac{1}{|C|} \sum_{c \in C} h(cg) - h(g) = \Big\langle h \Big| \tfrac{1}{|C|} \sum_{c \in C} P^n_{cg} - P^n_g \Big\rangle = \langle h \mid \nabla^* \tau_n \rangle = \langle \nabla h \mid \tau_n \rangle
\]
where $\tau_n$ is the transport plan obtained by taking the mass of $P^n_g$ at $g'$, split in $|C|$ masses, and take them (along a shortest path) to $g'c$ (for $c \in C$). Each transport takes at most $K:= \max_{c \in C} |c|_S$ steps. Notice that $g'C = Cg'$, so that this also splits and transports the mass uniformly to $Cg'$. Note that $\|\tau_n\|_{\ell^1} \leq K$ and that its $\ell^1$ norm on a ball of radius $k$ (denoted by $B_k$) tends to $0$ as $n \to \infty$ (it is bounded above by $K \| P^n_g \|_{\ell^1( B_{k+K})}$). On other hand, $\|\nabla h\|_{\ell^\infty(B_k^\mathsf{c})}$ tends to $0$ for $k$ large enough. Hence
\[
\begin{array}{r@{\,}l}
\big| \tfrac{1}{|C|} \sum_{c \in C} h(cg) - h(g) \big| 
&= |\langle \nabla h \mid \tau_n \rangle| \\
& \leq |\langle \nabla h \mid \tau_n \rangle_{B_k}| + |\langle \nabla h \mid \tau_n \rangle_{B_k^\mathsf{c}}|\\
& \leq \|\nabla h\|_{\ell^\infty} \|\tau_n\|_{\ell^1(B_k)} + \|\nabla h\|_{\ell^\infty(B_k^\mathsf{c})} \|\tau_n\|_{\ell^1}\\
& \leq \|\nabla h\|_{\ell^\infty} \|\tau_n\|_{\ell^1(B_k)} + \|\nabla h\|_{\ell^\infty(B_k^\mathsf{c})} K
\end{array}
\]
Note that the left-hand side does not depend on $n$ or $k$. By letting $n \to \infty$, the first term of the right hand-side tends to $0$. Then let $k \to \infty$, to see that the other term of the right-hand side tends to $0$. Hence $\sum_{c \in C} h(cg) = |C| h(g)$ for any finite conjugacy class $C$.

Pick some $g \in G$ and consider $a:= h(gs) - h(g)$ for some $s \in S$. Then $\sum_{c \in C} \big( h(cgs) - h(cg) \big) = |C|a$ for any finite conjugacy class $C$. In particular, for any finite conjugacy class there is a $c \in C$ so that $|h(cgs) - h(cg)| \geq |a|$. Since there are infinitely many finite conjugacy classes and $h$ has gradient in $C_0$, $a=0$. This shows that $h$ is constant.
\end{proof}
This lemma directly implies Proposition \ref{PasHarmRepC0}. We note that virtually nilpotent groups have infinitely many finite conjugacy classes, and therefore the previous lemma applies. Moreover, we mention that the conclusion of Proposition \ref{PasHarmRepC0} holds for 
\begin{itemize}
 \item any unitary representation of a nilpotent groups which does not contain the trivial representation, see \cite[Th\'eor\`eme 7 in \S{}8]{Gu}.
 \item any irreducible unitary representation of the Euclidean isometry groups, see \cite[Exemple 2 in \S{}9]{Gu}.
\end{itemize}
The case of polycyclic (and more generally solvable) groups is still open, even for $C_0$ representations.

\section{Proof of Theorem \ref{compression}}

\subsection{Proof of Claim (i) and (ii)}\label{secrecrel}

We start with a lemma.

\begin{Lem}\label{lapandgradient}
Let $X=(V,E)$ be a locally finite graph. Let $F : V \rightarrow \Hi$ be any map. 
\begin{enumerate}
	\item For any $x\in V$, the following identity holds:
		$$
		\mathcal{L}(\|F\|^{2})(x)=\|\nabla_{x}F\|_{\Hi}^{2}+2\Re\langle 						\mathcal{L}F(x),F(x)\rangle_{\Hi},
		$$
		where, by definition, 
		$$
		\|\nabla_{x}F\|^{2}=\|\nabla_{x}F\|_{\Hi}^{2}:=\frac{1}{\mathrm{deg}(x)}\sum_{y\sim x} \|F(x)-F(y)\|^{2}.
		$$
	\item Let $G$ be a group acting on $X$ by automorphisms and on $\Hi$ by affine isometries. Then, for any $G$-equivariant map $F : V \rightarrow \Hi$, the quantity $\|\nabla_{x}F\|$ is constant along each orbit of $x\in V$, that is,
	$$
	\|\nabla_{gx}F\|=\|\nabla_{x}F\|, \ \ \forall x\in V, g\in G. 
	$$
\end{enumerate}
\end{Lem}

\medskip

\noindent\textbf{Proof : } The proof of the first identity is a simple calculation. Since we have
$$
\|\nabla_{x}F\|^{2}=\| F(x)\|^{2}+\frac{1}{\mathrm{deg}(x)}\sum_{y\sim x} \| F(y)\|^{2}-\frac{2}{\mathrm{deg}(x)}\sum_{y\sim x}\Re \langle F(y),F(x)\rangle,
$$
and
$$
\langle\mathcal{L}F(x),F(x)\rangle = -\|F(x)\|^{2}+\frac{1}{\mathrm{deg}(x)}\sum_{y\sim x} \langle F(y),F(x)\rangle,
$$
it readily follows that
\begin{eqnarray*}
\|\nabla_{x}F\|^{2}+2\Re \langle \mathcal{L}F(x),F(x) \rangle
&=&
-\| F(x)\|^{2} + \frac{1}{\mathrm{deg}(x)}\sum_{y\sim x} \| F(y)\|^{2} \\
&=& 
\mathcal{L}\| F\|^{2}(x). 
\end{eqnarray*}

To prove the second statement, we recall that the equivariance implies
\begin{eqnarray*}
\| F(gx)-F(y)\|
&=&
\| \alpha(g)F(x)-F(y)\| \\
&=&
\| F(x)-\alpha(g^{-1})F(y)\| \\
&=&
\| F(x)-F(g^{-1}y)\| ,
\end{eqnarray*}
for all $g\in G$, $x,y\in V$, where $\alpha$ is the map corresponding to the $G$-action on $\Hi$. Hence, we obtain
\begin{eqnarray*}
\hspace*{3.05cm}
\|\nabla_{gx}F\|^{2}
&=&
\frac{1}{\mathrm{deg}(x)}\sum_{y\sim gx} \| F(gx)-F(y)\|^{2} \\
&=&
\frac{1}{\mathrm{deg}(x)}\sum_{y\sim gx} \| F(x)-F(g^{-1}y)\|^{2} \\
&=&
\|\nabla_{x}F\|^{2}. \hspace*{6.05cm} \square
\end{eqnarray*}


\begin{Lem}
Let $G$ be acting transitively on $\partial \tk$ and let us denote by $G_{x_{0}}$ the stabiliser of $x_{0}$. If the map $F : V \rightarrow \Hi$ satisfies $F(x_{0})=0$ for some point $x_{0}$ and if $F$ is $G$-equivariant with respect to some affine isometric action $\alpha=(\pi,b)$, where $b$ is zero on $G_{x_{0}}$, then both $\|F\|$ and $\langle\mathcal{L}F,F\rangle$ are radial.
\end{Lem}

\medskip

\noindent\textbf{Proof : } Let $x,y\in V$ be two vertices at distance $r$ from $x_{0}$. Since $G_{x_{0}}$ acts transitively on any sphere about $x_{0}$, we can find $h\in G_{x_{0}}$ so that $hx=y$. By the hypothesis on $b$, $\alpha$ acts by unitary operators when restricted to the stabiliser of $x_{0}$. Using the equivariance of $F$, it is straightforward that
$$
\|F(y)\|=\|\pi(h)F(x)\|=\|F(x)\|.
$$ 
Moreover, using the $G$-equivariance of $\mathcal{L}$ with respect to $\pi$, we obtain
\begin{eqnarray*}
\hspace*{2.8cm}
\langle \mathcal{L}F(y),F(y)\rangle 
&=&
\langle \mathcal{L}F(hx),F(hx)\rangle \\
&=&
\langle \pi(h)\mathcal{L}F(x),\pi(h)F(x) \rangle \\
&=&
\langle \mathcal{L}F(x),F(x) \rangle . \hspace*{4.3cm} \square
\end{eqnarray*}


\smallskip

For the rest of the section, we will assume that $G$, $F$ and $\alpha=(\pi,b)$ satisfy the hypothesis of Theorem \ref{compression}. Select one vertex on each sphere of radius $n$ and denote it by $x_{n}$. We can make two remarks. 
\begin{enumerate}
	\item Since $G$ acts transitively on $V$, then $\|\nabla_{x}F\|=\|\nabla_{x_{0}}F\|=\|F(x_{1})\|$.
	\item We recall that a non-compact closed subgroup of $\mathrm{Aut}(\tk)$ acts transitively on $\bordtk$ if and only if there exists $y\in V$ so that the stabiliser $G_{y}$ acts transitively on $\bordtk$ (see Proposition 10.1 in \cite{FTN}). If the group $G$ acts transitively on $V$, then the latter condition is also equivalent to having all the vertex stabilisers acting transitively on $\bordtk$. This means that under the hypotheses of Theorem \ref{compression}, the subgroup $G_{x_{0}}$ acts transitively on $\bordtk$.
\end{enumerate}

\medskip

Now, set $\varphi(n)=\|F(x_{n})\|^{2}$. For radial maps, the Laplace operator takes a very simple form :
$$  
\mathcal{L}\varphi(0)=\varphi(1)-\varphi(0),
$$
and
$$
\mathcal{L}\varphi(n)=\frac{q}{q+1}\varphi(n+1)-\varphi(n)+\frac{1}{q+1}\varphi(n-1),
$$
for any $n\geq 1$. Using Lemma \ref{lapandgradient} and setting $R_{F}(n)=2\Re\langle\mathcal{L}F(x_{n}),F(x_{n})\rangle$, we obtain the relations for all $n\geq 1$,
\begin{eqnarray}\label{recurrence}
\frac{q}{q+1}\varphi(n+1)-\varphi(n)+\frac{1}{q+1}\varphi(n-1)=\varphi(1)+R_{F}(n),
\end{eqnarray}
with initial conditions $\varphi(0)=0$ and $\varphi(1)=\|F(x_{1})\|^{2}$.\\
Let us find a general solution to this second order linear recurrence equation. Set
$$
\psi(n+1)=\varphi(n+1)-\varphi(n),
$$
for all $n\geq 0$. We can express the left-hand side of the relation (\ref{recurrence}) using $\psi$:
\begin{eqnarray*}
\frac{q}{q+1}\varphi(n+1)-\varphi(n)+\frac{1}{q+1}\varphi(n-1)
&=&
\frac{q}{q+1}(\varphi(n+1)-\varphi(n)) \\
& &
-\frac{1}{q+1}(\varphi(n)-\varphi(n-1)) \\
&=&
\frac{q}{q+1}\psi(n+1)-\frac{1}{q+1}\psi(n).
\end{eqnarray*}

For all $n\geq 1$, we obtain the new relation: 
$$
\psi(n+1)=\frac{1}{q}\psi(n)+\frac{(q+1)\|F(x_{1})\|^{2}}{q}+\frac{q+1}{q}R_{F}(n),
$$
with initial condition $\psi(1)=\|F(x_{1})\|^{2}$. By iterating this relation, we get
$$
\psi(n+1)=\frac{(q+1)\|F(x_{1})\|^{2}}{q-1}-\frac{2\|F(x_{1})\|^{2}}{(q-1)q^{n}}+(q+1)\sum_{j=1}^{n}\frac{R_{F}(n+1-j)}{q^{j}}.
$$
To proceed, we need a crucial negativity result.

\begin{Lem}\label{NegaLap} Under the hypotheses of Theorem \ref{compression}, we have 
$$
\Re\langle\mathcal{L}F(x),F(x)\rangle \leq 0,
$$
for all $x\in V$.
\end{Lem}

\medskip

Let us postpone the proof of this lemma to the next subsection and let us show how to finish the proof of Theorem \ref{compression}. Firstly, Lemma \ref{NegaLap} implies the following inequality
$$
\psi(n+1)\leq\frac{(q+1)\|F(x_{1})\|^{2}}{q-1}-\frac{2\|F(x_{1})\|^{2}}{(q-1)q^{n}},
$$
with equality if $F$ is harmonic. Replacing $\psi$ by $\varphi$, we get
$$
\varphi(n+1)\leq \varphi(n)+\frac{(q+1)\|F(x_{1})\|^{2}}{q-1}-\frac{2\|F(x_{1})\|^{2}}{(q-1)q^{n}},
$$ 
Iterating this inequality, we obtain the desired upper bound. Once again, if $F$ is harmonic, then the equality occurs, and we are done with the proof of the first two claims of Theorem \ref{compression}.\hfill $\square$\\

\smallskip

\subsection{Proof of Lemma \ref{NegaLap}}

The end of the present section is dedicated to the proof of Lemma \ref{NegaLap}. The main steps of the proof are described in the next lemma.

\begin{Lem} We assume that the group $G$ and the map $F$ satisfy the hypotheses of Theorem \ref{compression}. Also, we write $K$ for the compact subgroup $G_{x_{0}}$ and $dk$ for the normalised Haar measure on $G$ so that the subgroup $K$ has measure $1$. Fix $x\in V$ and let $g$ and $s$ in $G$ be such that $x=g^{-1}x_{0}$ and $sx_{0}$ is adjacent to $x_{0}$. 
	\begin{enumerate}
		\item\label{claim 1} For any vertex $y\in V$ adjacent to $x_{0}$, the map $F$ satisfies the following integral formula :
		$$
		\int_{K} F(g^{-1}ky) dk = \frac{1}{q+1}\sum_{u\sim x} F(u).
		$$
		\item\label{claim 2} Furthermore, if we denote by $\alpha=(\pi,b)$ the affine $G$-action on $\Hi$, then
		$$
		\mathcal{L}F(x)=\int_{K} \pi(g^{-1}k)b(s) dk.
		$$
		\item\label{claim 3} We also have
		$$
		\left\langle\mathcal{L}F(x), F(x)\right\rangle = -\left\langle P_{K}b(s),b(g) \right\rangle,
		$$
		where $P_{K}=\int_{K}\pi(k) dk$ is the orthogonal projection onto the space of $\pi(K)$-fixed vectors in $\Hi$.\\ 
		\item\label{claim 4} In particular, disintegrating $\pi$ and $b$ as direct integrals over some measure space $(Z,\nu)$, we obtain
		$$
		\langle\mathcal{L}F(x),F(x)\rangle = -\int_{Z}\langle P_{z}b_{z}(s),b_{z}(g)\rangle d\nu(z),
		$$
		where, for almost all $z\in Z$, $\pi_{z}$ is an irreducible unitary representation of $G$, $b_{z}$ is a $\pi_{z}$-cocycle, and 
		$$
		P_{z}=\int_{K}\pi_{z}(k) dk
		$$
		is the orthogonal projection onto the space of $\pi_{z}(K)$-invariant vectors.\\ 
		\item\label{claim 5} Finally, for any irreducible unitary representation $\sigma$ and for any cocycle $w\in Z^{1}(G,\sigma)$, we have
		$$
		\left\langle Pw(s), w(g)\right\rangle \geq 0,
		$$
		where $P=\int_{K}\sigma(k) dk$ denotes the orthogonal projection onto the space of $\sigma(K)$-invariant vectors.\\
	\end{enumerate}
\end{Lem} 

\medskip

Clearly, Lemma \ref{NegaLap} follows directly from the last two claims above.\\

\noindent\textbf{Proof of Claim \ref{claim 1}} : Fix $x\in V$. Since $G$ acts on $V$ transitively, there exists $g\in G$ so that $x=g^{-1}x_{0}$. Let $y_{1}, \dots, y_{q+1}$ be the $q+1$ neighbours of $x_{0}$. Using the fact that the action of $G_{x_{0}}$ on the sphere of radius $1$ centred at $x_{0}$ is transitive, we can find $h_{j}\in G_{x_{0}}$ so that $h_{j}y_{j}=y_{1}$, for $j=1,\dots, q+1$. We remark that the cosets given by $h_{j}^{-1}(G_{x_{0}}\cap G_{y_{1}})$ are all distinct and that the subgroup $G_{x_{0}}\cap G_{y_{1}}$ has index $q+1$ in $G_{x_{0}}$. Normalising the Haar measure on $G$ so that the compact subgroup $G_{x_{0}}$ has measure $1$, we obtain the following relation:  
\begin{eqnarray*}
\int_{G_{x_{0}}} F(g^{-1}ky_{1}) dk
&=&
\sum_{j=1}^{q+1} \int_{G_{x_{0}}\cap G_{y_{1}}} F(g^{-1}h_{j}^{-1}ky_{1}) dk \\
&=&
\sum_{j=1}^{q+1} F(g^{-1}y_{j}) \int_{G_{x_{0}}\cap G_{y_{1}}} dk \\
&=&
\frac{1}{q+1}\sum_{j=1}^{q+1}F({g^{-1}y_{j}})\\
&=&
\frac{1}{q+1}\sum_{u\sim x} F(u).
\end{eqnarray*}
In particular, we see that the integral on the left-hand side does not depend on the choice of the neighbour of $x_{0}$.\\

\smallskip
 
\noindent\textbf{Proof of Claim \ref{claim 2} : }   Again, by transitivity of the $G$-action, we can find $s\in G$ so that $sx_{0}=y_{1}$. It is straightforward to see that $b$, the translation part of the $G$-action on $\Hi$, factors through a $G$-equivariant map $G/K \rightarrow \Hi$ and that it coincides with $F$. Namely, $F(hx_{0})=b(h)$, for all $h\in G$. In particular, $F(x)=b(g^{-1})$. By the cocycle relation, we have
\begin{eqnarray*}
F(g^{-1}ky_{1})-F(g^{-1}kx_{0})
&=&
F(g^{-1}ksx_{0})-F(g^{-1}kx_{0})\\
&=&
b(g^{-1}ks)-b(g^{-1}k) \\
&=&
\pi(g^{-1}k)b(s)+b(g^{-1}k)-b(g^{-1}k) \\
&=&
\pi(g^{-1}k)b(s).
\end{eqnarray*} 

Hence, we obtain the desired integral formula for the Laplace operator applied to $F$:
\begin{eqnarray*}
\mathcal{L}F(x)
&=&
\left(\frac{1}{q+1}\sum_{y\sim x} F(y)\right)-F(x) \\
&=&
\left(\int_{K}F(g^{-1}ksx_{0}) dk\right)-F(g^{-1}x_{0}) \\
&=&
\int_{K}\left(F(g^{-1}ksx_{0})-F(g^{-1}x_{0}) \right)dk \\
&=&
\int_{K} \pi(g^{-1}k)b(s) dk.
\end{eqnarray*}\\

\smallskip

\noindent\textbf{Proof of Claim \ref{claim 3} : } Hence, applying Claim \ref{claim 2}, we directly get :
\begin{eqnarray*}
\left\langle\mathcal{L}F(x),F(x)\right\rangle
&=&
\left\langle\int_{K}\pi(g^{-1}k)b(s) dk, b(g^{-1}) \right\rangle \\
&=&
\int_{K}\left\langle \pi(g^{-1}k)b(s),b(g^{-1}\right)\rangle dk \\
&=&
\int_{K}\left\langle \pi(k)b(s),\pi(g)b(g^{-1})\right\rangle dk \\
&=&
-\int_{K}\left\langle \pi(k)b(s),b(g)\right\rangle dk \\
&=&
-\left\langle\int_{K}\pi(k)b(s)dk,b(g)\right\rangle \\
&=&
-\left\langle P_{K}b(s),b(g)\right\rangle.
\end{eqnarray*}\\

\smallskip

\noindent\textbf{Proof of Claim \ref{claim 4} : } Now, let us disintegrate $\pi$ and $b$ as direct integrals. We write
$$
\pi=\int_{Z}^{\oplus} \pi_{z} d\nu(z),
$$  
and 
$$
b=\int_{Z}^{\oplus} b_{z} d\nu(z),
$$
for some measure space $(Z,\nu)$. Recall that, for almost all $z\in Z$, $\pi_{z}$ is an irreducible unitary representation of $G$ and $b_{z}$ is a $\pi_{z}$-cocycle. Thus, we deduce from Claim \ref{claim 3} :
$$
\langle\mathcal{L}F(x),F(x)\rangle = -\int_{Z}\langle P_{z}b_{z}(s_{0}),b_{z}(g)\rangle d\nu(z).
$$\\

\smallskip

\noindent\textbf{Proof of Claim \ref{claim 5} : } As $G$ acts transitively on $\bordtk$, the couple $(G,K)$ forms a Gelfand pair (see chapter II, section 4 in \cite{FTN}) and we treat three cases. If the irreducible representation $\sigma$ is not spherical, then $P=0$ and the inner-product is $0$. If $\sigma=1_{G}$ is the trivial representation, then $w=0$ and the inner-product is again equal to $0$. Indeed, under the assumptions on $G$, we have that $H^{1}(G,1_{G})=Hom(G,\C)=0$ (see p. 5 of \cite{N})\footnote{We  recall the argument briefly for the reader's convenience. Let $\phi\in Hom(G,\C)$ be a continuous homomorphism. The group $G$ is generated by the vertex stabiliser $G_{x_{0}}$ and by the edge stabiliser $G_{\lbrack x_{0}, y_{1}\rbrack}$. Both subgroups being compact, their images under $\phi$ are compact subgroups of $\C$. Therefore, $\phi$ is identically $0$.}. Finally, if we suppose that $\sigma$ is spherical and non-trivial, then there exists a $\sigma(K)$-invariant vector $\eta$ such that 
$$
w(h)=(\sigma(h)-1)\eta, \ \forall h\in G. 
$$ 
Up to rescaling $w$, we can assume that $\eta$ has norm $1$. As the space of $\sigma(K)$-invariant vectors is one-dimensional, $P$ is a rank-one operator. So, we can write
$$
P\xi=\langle\xi,\eta\rangle\eta,
$$
for all $\xi\in\Hi$. This yields
\begin{eqnarray*}
\langle Pw(s),w(g)\rangle
&=&
\langle P(\sigma(s)-1)\eta, (\sigma(g)-1)\eta \rangle \\
&=&
\langle (\sigma(s)-1)\eta,\eta\rangle \langle\eta, (\sigma(g)-1)\eta \rangle \\
&=&
(\langle \sigma(s)\eta,\eta\rangle-1)(\langle\eta,\sigma(g)\eta\rangle -1)\\
&=&
\left(\phi(s)-1\right)\left(\phi(g^{-1})-1\right),
\end{eqnarray*}
where $\phi$ is the (normalised) positive-definite function of $\pi$ associated to $\eta$. In particular, $\phi$ is a spherical function, that is, a radial eigenfunction of the normalised adjacency operator on $\tk$ and $\phi(x_{0})=0$. This function being positive-definite and radial, it is therefore real-valued. By Cauchy-Schwarz, we conclude that $\phi(h)-1\leq 0$, for all $h\in G$ and therefore, the scalar product we began with is always positive or null. This ends the proof of the Claim \ref{claim 5}. \hfill $\square$\\

\smallskip

\subsection{Proof of Claim (iii)}\label{secequivharm}

The existence of an equivariant and harmonic map follows directly from Nebbia's result and from the orthogonal decomposition given in Proposition \ref{PropGuich}. Indeed, let $\sigma^{-}$ be the irreducible representation with non-trivial cohomology which appears in Theorem \ref{ThmNeb}. Fix two adjacent vertices $x_{0},x_{1}$, set $K=\textrm{stab}(x_{0})$ and let $b$ be a unbounded cocycle in $Z_{K}^{1}(G,\sigma^{-})$. Since Nebbia shows that $\sigma^{-}$ is a non-spherical representation, it implies that the subspace $\Hi^{\sigma^{-}(K)}$ is trivial. By Corollary \ref{CorNonSphericalHarmonic}, we get that $b$ is necessarily harmonic. Factorising $b$ through $G/K$, we obtain a $G$-equivariant map $F$ which is harmonic and non-constant. This finishes the proof of Theorem \ref{compression}. \hfill $\square$ 

\medskip

In fact, since $\sigma^{-}$ is the unique irreducible representation of $G$ which has non-vanishing cohomology and since $H^{1}(G,\sigma^{-})$ has complex dimension $1$, the cocycle $b$ appearing in the proof here above is unique up to scalar multiplication. Therefore, Theorem \ref{compression} gives an alternative description of the non-trivial cocycle appearing at the end of Nebbia's paper.

\bigskip

\section{The Haagerup cocycle}

The goal of the present section is to decompose the so-called \textbf{Haagerup cocycle} and to give an alternative description of the proper harmonic cocycle appearing in Theorem \ref{compression}. We first need to introduce some notations (for all this, see p. 90 of \cite{BHV} and Chapter 1 of \cite{Woe}).

\subsection{Decomposition of the Haagerup cocycle}\label{decHaag}

 Let $X=(V,\mathbb{E})$ be a locally finite graph, where $\mathbb{E}$ denotes the set of \textbf{oriented} edges. Each edge $e\in\mathbb{E}$ has a source $s(e)=e_{-}\in V$ and a range $r(e)=e_{+}\in V$. There is an obvious fixed-point free involution $e\mapsto \overline{e}$ on $\mathbb{E}$ with $s(\overline{e})=r(e)$ and $r(\overline{e})=s(e)$, for all $e\in\mathbb{E}$. The set of all pairs $\{e,\overline{e}\}$ is the set of \textbf{geometric edges} of the graph $X$. We denote by $\ell_{\mathrm{alt}}^{2}(\mathbb{E})$ the real Hilbert space of those maps $\xi : \mathbb{E} \rightarrow \R$ satisfying $\xi(\overline{e})=-\xi(e)$ and such that $\sum_{e\in\mathbb{E}}|\xi(e)|^{2}<\infty$. This vector space is endowed with the inner product
$$
\langle \xi,\eta \rangle=\frac{1}{2}\sum_{e\in\mathbb{E}}\xi(e)\eta(e).
$$
Let $\ell^{2}(V,\mathrm{deg})$ be the Hilbert space of square-summable functions on $V$ endowed with the inner product
$$
(f,g)=\sum_{x\in V}f(x)g(x)\mathrm{deg}(x).
$$
Now we can define two operators connecting these spaces. Let $\nabla : \ell^{2}(V,\mathrm{deg})\rightarrow \ell_{\mathrm{alt}}^{2}(\mathbb{E})$ be the \textbf{gradient}, defined by
$$
\left(\nabla f\right)(e)=f(e_{+})-f(e_{-}).
$$
It is straightforward from the definition that $\left(\nabla f\right)(\overline{e})=-\left(\nabla f\right)(e)$, for all $e\in\mathbb{E}$. We also define the \textbf{divergence} $\nabla^{\ast} : \ell_{\mathrm{alt}}^{2}(\mathbb{E})\rightarrow \ell^{2}(V,\mathrm{deg})$ as the adjoint of $\nabla$, that is,
$$
\left(\nabla^{\ast}u\right)(x)=\frac{1}{\mathrm{deg}(x)}\sum_{y\sim x}u(y,x),
$$
for all $u\in\ell_{\mathrm{alt}}^{2}(\mathbb{E})$ and $x\in V$. The link between the Laplace operator $\mathcal{L}$ and the operators $\nabla$ and $\nabla^{*}$ is given by the following formula:
$$
\mathcal{L}=-\nabla^{*}\nabla.
$$
\begin{Lem}(Poincaré lemma on trees)
Let $\mathcal{T}=(V,\mathbb{E})$ be a tree and fix a vertex $x_{0}$. For any map $\xi : \mathbb{E}\rightarrow \R$ such that $\xi(\overline{e})=-\xi(e)$, for all $e\in\mathbb{E}$, there is a unique function $\tilde{\xi} : V\rightarrow \R$ such that $\nabla\tilde{\xi}=\xi$ and $\tilde{\xi}(x_{0})=0$.
\end{Lem}

\medskip

\noindent {\bf Proof : } Set $\tilde{\xi}(x_{0})=0$. Let $n\geq 1$. Let $x_{n}$ be a vertex at distance $n$ from $x_{0}$. Let $(x_{j})_{j=0}^{n}$ be the chain of vertices forming the unique geodesic path linking $x_{0}$ to $x_{n}$ in $X$. Set 
$$
\tilde{\xi}(x_{n})=\sum_{j=0}^{n-1}\xi(x_{j},x_{j+1}).
$$
It is easy to see that $\nabla\tilde{\xi}=\xi$. \hfill $\square$\\

\smallskip

Recall that a function on the vertices $\eta : V \rightarrow \R$ is harmonic if $\mathcal{L}\eta=0$. In the case where $X=\mathcal{T}$ is a tree, using the previous lemma, a map $\xi\in\ell^{2}(\mathbb{E})$ belongs to $\ker\nabla^{\ast}$ if and only if $\tilde{\xi}$ is harmonic. Therefore, it is natural to say that a map on the edges $\xi\in\ell_{\mathrm{alt}}^{2}(\mathbb{E})$ is \textbf{harmonic} if $\xi\in\ker\nabla^{\ast}$. This remark suggests the following orthogonal decomposition: 
$$
\ell_{\mathrm{alt}}^{2}(\mathbb{E})=\ker\nabla^{\ast}\oplus \overline{\mathrm{im}\nabla},
$$ 
since $(\ker\nabla^{\ast})^{\perp}=\overline{\mathrm{im}\nabla}$. We can give a more precise description of the projection onto $\overline{\mathrm{im}\nabla}$. Recall first that, if the graph is non amenable, then $\mathrm{im}\nabla$ is closed and the Laplace operator $\mathcal{L}$ is invertible. Furthermore, its inverse is the operator $-G$, with $G$ being the \textbf{Green kernel} of $X$ (see p.14 of \cite{Woe} for the definition). Let $\xi=h+\nabla k$, with $h\in\ker\nabla^{\ast}$ and $k\in\ell^{2}(V,\mathrm{deg})$. Then, $\nabla^{\ast}\xi=\nabla^{\ast}\nabla k=-\mathcal{L}k$. Hence, we obtain $k=G\nabla^{\ast}\xi$. This forces to define $Q$, the orthogonal projection onto $\mathrm{im}(\nabla)$, by 
$$
Q(\xi)=\nabla G\nabla^{\ast}\xi.
$$
We can deduce that the ``harmonic part'' of an element $\xi\in\ell_{\mathrm{alt}}^{2}(\mathbb{E})$ is given by $(1-Q)\xi$. Now, we can introduce the Haagerup cocycle and study its decomposition along the two subspaces $\ker\nabla^{*}$ and $\mathrm{im}\nabla$.\\
Let $X=(V,\mathbb{E})$ be a graph. For two vertices $x,y\in V$, we define the \textbf{signed characteristic function} of the geodesic $\lbrack x,y\rbrack$ by
$$
\chi_{x\rightarrow y}(e)=\left\{
	\begin{array}{lll}
		1, & \textrm{if $e$ is on $\lbrack x,y\rbrack$ and $e$ points from $x$ to $y$}, \\
		-1, & \textrm{if $e$ is on $\lbrack x,y\rbrack$ and $e$ points from $y$ to $x$}, \\
		0, & \textrm{otherwise},
	\end{array}\right.
$$
In the case where $X=\mathcal{T}$ is a tree, then a simple calculation shows that
$$
\|\chi_{x\rightarrow y}\|_{\ell_{\mathrm{alt}}^{2}(\mathbb{E})}= \sqrt{d(x,y)}.
$$ 
Let $G$ be a closed subgroup of $\mathrm{Aut}(\mathcal{T})$. Let $\pi$ be the orthogonal representation of $G$ on $\ell_{\mathrm{alt}}^{2}(\mathbb{E})$ induced by the action of $G$ on $\mathcal{T}=(V,\mathbb{E})$. Let $x_{0}\in V$ be fixed.
The \textbf{Haagerup cocycle} is defined by $b : G \rightarrow \ell_{\mathrm{alt}}^{2}(\mathbb{E})$ with
$$ 
b(g)=\chi_{x_{0}\rightarrow gx_{0}},
$$
It is easy to check that $b$ satisfies the cocycle relation with respect to $\pi$. By the previous observation, we have the following identity:
$$
\|b(g)\|_{\ell_{\mathrm{alt}}^{2}(\mathbb{E})}=\sqrt{d(x_{0},gx_{0})}.
$$
This proves that $b$ is a proper cocycle and that $G$ has the Haagerup property. Since $\ell_{\mathrm{alt}}^{2}(\mathbb{E})$ can be decomposed into an orthogonal sum of two (closed) $G$-invariant subspaces, then we can conclude that the representation $\pi$ is reducible. However, we will show in the sequel that the projection of the cocycle $b$ onto $\ker\nabla^{\ast}$ is still proper. To do so, we will prove that $\|Q\chi_{x\rightarrow y}\|_{\ell_{\mathrm{alt}}^{2}(\mathbb{E})}$ is bounded, independently of $x$ and $y$.\\ 
Here is a useful lemma allowing us to estimate the operator norm of the Green kernel.

\begin{Lem}\label{normG} Let $X=(V,E)$ be a graph and let $P$ be the normalised adjacency operator acting on $\ell^{2}(V,\mathrm{deg})$, namely, the operator whose matrix coefficients $p(x,y)$ are
$$
p(x,y)=\left\{
	\begin{array}{ll}
		\frac{1}{\mathrm{deg}(x)}, & \textrm{if $x\sim y$}, \\
		0, & \textrm{otherwise}.	
	\end{array}\right.
$$
\begin{enumerate}
	\item[(i)]If $\|P\|<1$, then the series $\sum_{n\geq 0}P^{n}$ defines a bounded operator and we have the equalities
$$
\mathcal{L}^{-1}=-\sum_{n\geq 0}P^{n}=-G,
$$
where $G$ is the Green kernel.\\
In particular, $\|G\|\leq\frac{1}{1-\|P\|}$.
	\item[(ii)] (Theorem (11.1), \cite{Woe}) If $X$ is a graph with all vertices of valency bounded by $q+1$, then $\|P\|\geq \frac{2\sqrt{q}}{q+1}$, with equality if $X=\mathcal{T}_{q+1}$.
\end{enumerate}
\end{Lem}

\medskip

For $x,y\in V$, let us estimate the norm of the harmonic part of $\chi_{x\rightarrow y}$. Firstly, it is a general fact that, for any $\xi\in\ell_{\mathrm{alt}}^{2}(\mathbb{E})$, we have
$$
\|(1-Q)\xi\|\leq \|\xi\|,
$$
and
$$
\|(1-Q)\xi\|^{2}=\|\xi\|^{2}-\| Q(\xi)\|^{2}.
$$
Now, let us compute $Q\chi_{x\rightarrow y}$. It is easy to check that 
$$
\nabla^{\ast}(\chi_{x\rightarrow y})=\frac{\delta_{y}}{\mathrm{deg}(y)}-\frac{\delta_{x}}{\mathrm{deg}(x)}.
$$
We get
\begin{eqnarray*}
\|\nabla^{\ast}\chi_{x\rightarrow y}\|_{\ell^{2}(V,\mathrm{deg})}^{2}
&=&
\Big\|\frac{\delta_{y}}{\mathrm{deg}(y)}-\frac{\delta_{x}}{\mathrm{deg}(x)} \Big\|_{\ell^{2}(V,\mathrm{deg})}^{2} \\
&=&
\frac{1}{\mathrm{deg}(y)}+\frac{1}{\mathrm{deg}(x)}.
\end{eqnarray*}
Using Lemma \ref{normG}, we obtain
\begin{eqnarray*}
\|(1-Q)\chi_{x\rightarrow y}\|^{2}
&=& 
\|\chi_{x\rightarrow y}\|^{2}-\| Q(\chi_{x\rightarrow y})\|^{2} \\
&\geq &
d(x,y)- \|\nabla\|^{2}\|G\|^{2}\|\nabla^{\ast}\chi_{x\rightarrow y}\|^{2} \\
&\geq &
d(x,y)-\frac{2}{(1-\| P\|)^{2}}\left(\frac{1}{\mathrm{deg}(x)}+\frac{1}{\mathrm{deg}(y)}\right) \\
&\geq &
d(x,y)-\frac{4}{(1-\| P\|)^{2}},
\end{eqnarray*}
since $\|\nabla\|=\|\nabla^{\ast}\|=\sqrt{2}$.\\
In particular, the orthogonal projection of the cocycle $b$ onto $\ker\nabla^{\ast}$ is still proper and its compression exponent is $\frac{1}{2}$. Indeed, for any $g\in \mathrm{Aut}(\mathcal{T})$, we have
$$
\|(1-Q)b(g)\|^{2} \geq d(x_{0},gx_{0})-\frac{4}{(1-\|P\|)^{2}}.
$$
If $\mathcal{T}=\mathcal{T}_{q+1}$ is the homogeneous $(q+1)$-regular tree, then it is possible to compute $Q(\chi_{x\rightarrow y})$ explicitly. We will prove the following lemma.

\begin{Lem}\label{projection}
Let $q\geq 2$ be an integer and let $\mathcal{T}_{q+1}$ be the homogeneous $(q+1)$-regular tree. Then, for any $x,y\in V$, we have:
$$
\hspace*{0.5ex} \text{(i)} \hspace*{0.5ex} |(Q\chi_{x\rightarrow y})(e)|=\left\{
	\begin{array}{ll}
		\displaystyle \frac{q^{-d(y,e)}+q^{-d(x,e)}}{q+1}, & \textrm{if $e$ is on the geodesic $\lbrack x,y\rbrack$}, \\
		\displaystyle \left|\frac{q^{-d(y,e)}-q^{-d(x,e)}}{q+1}\right|, & \textrm{otherwise,}
	\end{array}\right. 
$$
\begin{enumerate}
\item[(ii)] Moreover,
$$
\|Q(\chi_{x\rightarrow y})\|_{\ell_{\mathrm{alt}}^{2}(\mathbb{E})}^{2}=\frac{2q}{q^{2}-1}(1-q^{-d(x,y)}).
$$
\end{enumerate}
\end{Lem}

\medskip

We immediately obtain :

\smallskip

\begin{Cor}\label{ProjectionCompression}
Let $G$ be a subgroup of $\mathrm{Aut}(\tk)$ satisfying the assumptions of Theorem \ref{compression}. Fix a basepoint $x_{0}\in V$, set $|g|:=d(x_{0},gx_{0})$, for $g\in G$, and write $\tilde{b}$ for the projection of the Haagerup cocycle onto the closed invariant subspace $\ker \nabla^{*}$, that is, $\tilde{b}(g)=(1-Q)\chi_{x_{0}\rightarrow gx_{0}}$. Then, the cocycle $\tilde{b}$ is proper and it satisfies the following estimate
	$$
	\|\tilde{b}(g)\|^{2}= |g|+\frac{2q}{q^{2}-1}(q^{-|g|}-1).
	$$
In particular, the cocycle $\tilde{b}$ attains the upper bound of Theorem \ref{compression}.(i).
\end{Cor}

\medskip

Before proving Lemma \ref{projection}, we recall that the Green kernel $G$ takes a particularly simple form on $\mathcal{T}_{q+1}$.

\begin{Lem}(Lemma (1.23), \cite{Woe})
Let $q\geq 2$ and let $G$ be the Green kernel defined on the homogeneous $(q+1)$-regular tree $\mathcal{T}_{q+1}$. If we denote by $(G(x,y))_{x,y\in V}$ the associated matrix of $G$, then
$$
G(x,y)=\frac{q^{1-d(x,y)}}{q-1},
$$
for all $x,y\in V$.
\end{Lem}

\medskip

\noindent {\bf Proof of Lemma \ref{projection} : } Since $\nabla^{*}\chi_{x\rightarrow y}=\frac{1}{q+1}(\delta_{y}-\delta_{x})$, we need to compute $\nabla G \delta_{x}$. Let $e\in\mathbb{E}$ be an oriented edge. Using the description of the Green kernel, we get
\begin{eqnarray*}
(\nabla G \delta_{x})(e)
&=&
(G\delta_{x})(e_{+})-(G\delta_{x})(e_{-}) \\
&=&
G(e_{+},x)-G(e_{-},x) \\
&=&
\frac{q}{q-1}\left(q^{-d(e_{+},x)}-q^{-d(e_{-},x)}\right).
\end{eqnarray*}
Clearly, $|d(e_{+},x)-d(e_{-},x)|=1$. Setting $d(x,e)=\min\{d(x,e_{-}),d(x,e_{+})\}$ (this is simply the natural distance between $x$ and the geometric edge associated with $e$ in the geometric realisation of $\mathcal{T}_{q+1}$), we immediately obtain:
$$
\left(\nabla G\delta_{x}\right)(e) = \left\{
		\begin{array}{ll}
			-q^{-d(x,e)}, & \textrm{if $d(x,e)=d(x,e_{-})$}, \\
			q^{-d(x,e)}, & \textrm{if $d(x,e)=d(x,e_{+})$}		
		\end{array} \right.
$$
It is easy to see that $\left(\nabla G\delta_{x}\right)(e)$ and $\left(\nabla G\delta_{y}\right)(e)$ have the same sign if and only if $x$ and $y$ belong to the same connected component of \footnote{Here, $mid(e)$ denotes the median point of $e$ in the geometric realisation of $\mathcal{T}_{q+1}$.} $\mathcal{T}_{q+1}\setminus\{mid(e)\}$, which happens exactly when $e$ does not lie on the geodesic $\lbrack x, y\rbrack$. This shows the first claim.\\
From the first claim, we deduce:
\begin{eqnarray*}
\|Q\chi_{x\rightarrow y}\|_{\ell_{\mathrm{alt}}^{2}(\mathbb{E})}^{2} 
&=&
\frac{1}{2}\sum_{e\in\mathbb{E}} |(Q\chi_{x\rightarrow y})(e)|^{2} \\
&=&
\sum_{e\in E} |(Q\chi_{x\rightarrow y})(e)|^{2} \\
&=&
\frac{1}{(q+1)^{2}}\sum_{e\in E} |\left(\nabla G\delta_{y}-\nabla G\delta_{x}\right)(e)|^{2}.
\end{eqnarray*}
To compute the last sum, we will decompose the set of geometric edges. First of all, let $m=d(x,y)$ and let $\{z_{j}\}_{j=0}^{m}$ be the set of vertices describing the geodesic $\lbrack x,y\rbrack$, with $z_{0}=x$ and $z_{m}=y$. Let $T_{0}$ be the subgraph which is induced on the connected component of $\mathcal{T}_{q+1}\setminus\{z_{1}\}$ containing $x$. For $1\leq j\leq m-1$, let $T_{j}$ be the subgraph which is induced on the connected component of $\mathcal{T}_{q+1}\setminus\{z_{j-1},z_{j+1}\}$ containing $z_{j}$. Finally, let $T_{m}$ be the subgraph which is induced on the connected component of $\mathcal{T}_{q+1}\setminus\{z_{m-1}\}$ containing $y$. We remark that for all $j$, the graph $T_{j}$ is a subtree of $\mathcal{T}_{q+1}$ with root $z_{j}$. With these notations, a geometric edge $e$ belongs either to one of the $T_{j}$, for some $j$, or $e$ lies on $\lbrack x,y\rbrack$. Thus, we get:
\begin{eqnarray*}
\|Q\chi_{x\rightarrow y}\|_{\ell_{\mathrm{alt}}^{2}(\mathbb{E})}^{2} 
&=&
\frac{1}{(q+1)^{2}}\left(\left(\sum_{e\in \lbrack x,y\rbrack} \left|q^{-d(y,e)}+q^{-d(x,e)}\right|^{2}\right) \right. \\
& &
\left. +\left(\sum_{j=0}^{m}\sum_{e\in T_{j}} \left|q^{-d(y,e)}-q^{-d(x,e)}\right|^{2}\right)\right).
\end{eqnarray*}

To compute the first sum, let us denote by $e_{j}$ the edge $(z_{j},z_{j+1})$. Therefore, we have
\begin{eqnarray*}
\sum_{e\in \lbrack x,y\rbrack} \left|q^{-d(y,e)}+q^{-d(x,e)}\right|^{2}
&=&
\sum_{j=0}^{m-1} \left|q^{-d(y,e_{j})}+q^{-d(x,e_{j})}\right|^{2} \\
&=&
\sum_{j=0}^{m-1} \left|q^{-(m-j-1)}+q^{-j}\right|^{2} \\
&=&
\frac{2m}{q^{m-1}}+2\frac{1-q^{-2m}}{1-q^{-2}}.
\end{eqnarray*}
 
Secondly, let us compute the sum over the edges belonging to the subtree $T_{0}$. For any edge $e$ in $T_{0}$, we notice that $d(y,e)=m+d(x,e)$. Since the number of edges in $T_{0}$ which are at distance $k$ to $x$ is equal to $q^{k+1}$, for $k\geq 0$, we have:
\begin{eqnarray*}
\sum_{e\in T_{0}} \left|q^{-d(y,e)}-q^{-d(x,e)}\right|^{2}
&=&
\sum_{k\geq 0} \sum_{e\in T_{0} : \atop d(x,e)=k} |q^{-k}-q^{-k-m}|^{2} \\
&=&
\left(1-q^{-m}\right)^{2}\sum_{k\geq 0} q^{-k+1} \\
&=&
\left(1-q^{-m}\right)^{2}\frac{q}{1-q^{-1}}.
\end{eqnarray*}

By symmetry, the same is true for the sum over $T_{m}$. That is:
$$
\sum_{e\in T_{m}} \left|q^{-d(y,e)}-q^{-d(x,e)}\right|^{2}=\left(1-q^{-m}\right)^{2}\frac{q}{1-q^{-1}}.
$$
Finally, we need to compute the sum over the edges belonging to $T_{j}$ for $1\leq j\leq m-1$. For any edge $e$ in $T_{j}$, we notice that $d(x,e)=d(x,z_{j})+d(z_{j},e)=j+d(z_{j},e)$ and $d(y,e)=d(y,z_{j})+d(z_{j},e)=m-j+d(z_{j},e)$. Since the number of edges in $T_{j}$ which are at distance $k$ to $z_{j}$ is equal to $(q-1)q^{k}$, for $k\geq 0$, we have:
\begin{eqnarray*}
\sum_{e\in T_{j}} \left|q^{-d(y,e)}-q^{-d(x,e)}\right|^{2}
&=&
\sum_{k\geq 0} \sum_{e\in T_{j} : \atop d(z_{j},e)=k} |q^{-m+j-k}-q^{-j-k}|^{2} \\
&=&
\sum_{k\geq 0} \sum_{e\in T_{j} : \atop d(z_{j},e)=k} q^{-2k}|q^{-m+j}-q^{-j}|^{2} \\
&=&
\left(q^{-m+j}-q^{-j}\right)^{2}\sum_{k\geq 0} (q-1)q^{k}q^{-2k} \\
&=&
q\left(q^{-2(m-j)}+q^{-2j}-2q^{-m}\right).
\end{eqnarray*}

We can compute the sum over all the $T_{j}$, for $1\leq j\leq m-1$:

\begin{eqnarray*}
\sum_{j=1}^{m-1}\sum_{e\in T_{j}} \left|q^{-d(y,e)}-q^{-d(x,e)}\right|^{2}
&=&
\sum_{j=1}^{m-1}q\left(q^{-2(m-j)}+q^{-2j}-2q^{-m}\right) \\
&=&
-\frac{2(m-1)}{q^{m-1}}+2\frac{q^{-1}-q^{1-2m}}{1-q^{-2}} \\
\end{eqnarray*}

Since we have

\begin{eqnarray*}
\frac{1}{2}\sum_{e\in E} |\left(\nabla G\right)\left(\delta_{y}-\delta_{x}\right)(e)|^{2}
&=&
\frac{1}{q^{m-1}}+\frac{1-q^{-2m}+q^{-1}-q^{1-2m}}{1-q^{-2}}+\frac{q(1-q^{-m})^{2}}{1-q^{-1}} \\
&=&
\frac{(q^{m}-1)(q+1)}{(q-1)q^{m-1}} ,
\end{eqnarray*}

we deduce finally that
\begin{eqnarray*}
\| Q\chi_{x\rightarrow y}\|^{2}
&=&
\frac{1}{(q+1)^{2}}\sum_{e\in E}|(\nabla G)(\delta_{y}-\delta_{x})(e)|^{2} \\
&=&
\frac{2q}{q^{2}-1}\left(1-q^{-m}\right),
\end{eqnarray*}
which proves the second claim.\hfill $\square$\\

\smallskip

\subsection{Virtual coboundaries}\label{VirtualCoboundaries}

The aim of the current subsection is to describe the Haagerup cocycle and its projection onto $\ker\nabla^{\ast}$ as virtual coboundaries, which is a classical trick to produce unbounded cocycles. In order to speak of virtual coboundary, the Hilbert space $\Hi$ has to lie in a vector space $W$ and the unitary representation $\pi$ has to extend to a linear group action on $W$. A cocycle is then a virtual coboundary if $b(g) = \pi(g) x - x$ for some $x \in W \setminus \Hi$.\\
The situation which concerns us for the present paper is when $G$ acts transitively on the countable set of edges $\mathbb{E}$ and $\pi$ is the natural representation on $\ell_{\textrm{alt}}^{2}(\mathbb{E})$. The natural choice for the space $W$ is simply the whole set of alternate functions on $\mathbb{E}$, \emph{i.e.} $W = \{ f: \mathbb{E} \rightarrow \R \  : \  f(\overline{e})=-f(e)\}$.\\
Let us first describe the Haagerup cocycle as a virtual coboundary of the form $\pi(g)f-f$ (see also \cite[\S{}3]{FV} for another possible choice). A simple computations shows that $f$ may be defined as follows:
\[
f(e) = \left\{\begin{array}{ll}
+1/2 & \text{if following } e \text{ increases the distance to } x_0\\
-1/2 & \text{if following } e \text{ decreases the distance to } x_0\\
\end{array} \right.
\]
Note that the divergence of this function is $\tfrac{1}{2}$ at $x_0$ and $\tfrac{q-1}{2(q+1)}$ at every other vertex. We remark that for each $g\in G$, the function $b(g) = \pi(g) f- f\in\ell_{\textrm{alt}}^{2}(\mathbb{E})$ has a non-trivial divergence only at the vertices $x_0$ and $gx_0$. By the previous section, the Haagerup cocycle $b$ is linked to the optimal harmonic cocycle $\tilde{b}$ by the following relation : $b(g)=\tilde{b}(g)+b'(g)$, where $b'(g)=Q\chi_{x_0 \to g x_0}\in W$. The computations of \S{}\ref{decHaag} show that $b'$ is a coboundary: $b'(g) = \pi(g)f' -f'$ where $f' = \nabla G \delta_{x_0}$. It is then apparent that the optimal cocycle $\tilde{b}$ is also a virtual coboundary: $\tilde{b}(g)= \pi(g) \tilde{f} - \tilde{f}$ with $\tilde{f} = f+f'$. Note that the function $\tilde{f}$ has divergence $\tfrac{q-1}{2(q+1)}$ at every vertex, so that $\tilde{b}(g) = \pi(g) \tilde{f}- \tilde{f}$ lies indeed in the kernel of the divergence.
Since we are on a tree, all these alternated functions on the edges can be integrated as functions on the vertices. For example, $f$ is $\tfrac{1}{2}$ times the gradient of the function $x \mapsto d(x,x_0)$. Note that $\tilde{f}$ being of constant divergence and spherical (both $f$ and $f'$ are spherical), its integral will have constant Laplacian and be spherical. In particular, this means it will satisfy the recurrence relation described in \S{}\ref{secrecrel}.

\section{Conditionally negative type functions on $G$}

In this last section, we exploit Theorem \ref{compression} to give a comprehensive description of pure elements of $\mathrm{CL}(G)$, the convex cone of conditionally negative type functions on $G$. In order to state the result, we recall a few facts about negative type kernels and negative type functions. We refer to the Appendix C of \cite{BHV} for more details.

\smallskip

\subsection{Definitions and application of Theorem \ref{compression}}\label{secdeccone}

A kernel $\Psi : W \times W \rightarrow \R$ on a set $W$ is said to be \textbf{conditionally of negative type} if it satisfies the following properties:

\begin{enumerate}
	\item[(i)] $\Psi(x,x)=0$, for all $x\in W$; \\
	\item[(ii)] $\Psi(x,y)=\Psi(y,x)$, for all $x,y\in W$; \\
	\item[(iii)] For any $x_{1}, \dots, x_{n}\in W$ and for any $\alpha_{1},\dots ,\alpha_{n}\in \R$ satisfying\\ $\sum_{i=1}^{n}\alpha_{i}=0$, we have
	$$
	\sum_{i=1}^{n}\sum_{j=1}^{n}\alpha_{i}\alpha_{j}\Psi(x_{i},x_{j}) \leq 0.
	$$
\end{enumerate} 

Examples of such kernels are given by maps of the form $\Psi(x,y)=\|f(x)-f(y)\|^{2}$, where $f : W \rightarrow \Hi$ is any map with image in a Hilbert space. A continuous function $\psi : G \rightarrow \R$ on a topological group $G$ is said to be \textbf{conditionally of negative type} if the kernel on $G$ defined by $(g,h) \mapsto \psi(g^{-1}h)$ is conditionally of negative type. It is easy to see that if $b$ is a continuous cocycle for some unitary representation $\pi$ of $G$, then the function $g \mapsto \|b(g)\|^{2}$ is conditionally of negative type. This example is essentially universal. Namely, given $\psi$ a function conditionally of negative type on a group $G$, by the GNS construction, there exist $\pi_{\psi}$ a unitary representation of $G$ and a cocycle $b_{\psi}\in Z^{1}(G,\pi_{\psi})$ satisfying $\psi(g)=\|b_{\psi}(g)\|^{2}$, for all $g\in G$. It is well-known that, for a group $G$, the set of such functions forms a convex, positive cone, denoted by $\mathrm{CL}(G)$. We say that a function $\psi\in \mathrm{CL}(G)$ is \textbf{pure} if it lies on an extremal ray of $\mathrm{CL}(G)$. We recall the following:

\begin{Thm}\label{LSV} (Théorème 1, \cite{LSV}) Let $G$ be a topological group.
\begin{enumerate}
	\item[(i)] Let $\psi$ be a function conditionally of negative type on $G$ and let $(\pi_{\psi},\Hi_{\psi},b_{\psi})$ be its associated GNS triple. If $\psi$ is pure, then the orthogonal representation $\pi_{\psi}$ is irreducible.\\
	\item[(ii)] Let $\pi$ be an irreducible orthogonal representation and let $b$ be any 1-cocycle for the representation $\pi$. Then, the function of negative type $\psi$ associated with $b$ is pure.
\end{enumerate}
\end{Thm}

\medskip

This yields:  

\begin{Cor}\label{UniqueNegTypeFun}
\begin{enumerate}
	\item[(i)] The kernel defined on the set of vertices of $\tk$ by
$$
\Psi : (x,y)\mapsto d(x,y)-\frac{2q}{q^{2}-1}+\frac{2q}{q^{2}-1}\cdot q^{-d(x,y)}
$$
is conditionally of negative type.\\
	\item[(ii)] Let $G$ be as in Theorem \ref{compression} and choose a basepoint $x_{0}\in V$. For $g\in G$, we set $|g|:=d(x_{0},gx_{0})$. Then, the function on $G$ defined by
$$
g \mapsto \Psi(gx_{0},x_{0})=|g|+\frac{2q}{q^{2}-1}(q^{-|g|}-1)
$$ 
is the unique (up to multiplication by a positive scalar) pure negative type function in $\mathrm{CL}(G)$ which is unbounded on $G$ and identically $0$ on $G_{x_{0}}$.\\
	\item[(iii)] The cocycle $\tilde{b}$ of Corollary \ref{ProjectionCompression} coincides with the harmonic equivariant map $F$ of Theorem \ref{compression}.(iii).
\end{enumerate}
\end{Cor}

\medskip

\noindent {\bf Proof : } To prove the first claim, we simply remark that 
$$
\Psi(x,y)=\|F(x)-F(y)\|^{2},
$$ 
where the map $F$ is as in Theorem \ref{compression}.(iii) and is normalised so that 
$$
\|F(x_{1})\|^{2}=\frac{q-1}{q+1}.
$$ 
To prove claim (ii), we need to see that the function $g \mapsto \Psi(gx_{0},x_{0})$ is pure. Let $\sigma^{-}$ be the unique irreducible unitary representation which admits an unbounded cocycle $h\in Z_{K}^{1}(G,\sigma^{-})$. Since this representation is non-spherical, we deduce that $h$ is harmonic, by Corollary \ref{CorNonSphericalHarmonic}. Therefore, the function conditionally of negative type associated with $h$ attains the upper bound of Theorem \ref{compression}.(i), which coincides with $\Psi(gx_{0},x_{0})$, up to a multiplicative constant. As a consequence of Theorem \ref{LSV}, this function is pure. By Nebbia's result, the space $Z_{K}^{1}(G,\sigma^{-})$ has one dimension and the choice of a cocycle $h$ here above is unique, up to a scalar multiplication. This proves the uniqueness of the function of negative type. Finally, by claim (ii) and Theorem \ref{LSV}, the equivariant maps $F$ of Theorem \ref{compression} and $\tilde{b}$ are cocycles relatively to the representation $\sigma^{-}$, and therefore, they must coincide. \hfill $\square$\\

\smallskip

\subsection{More examples and classification of pure conditionally negative type functions}\label{SubsectionExamplesAndCNTFun}

We start by giving an interesting family of examples of kernels conditionally of negative type on trees. Let $\mathcal{T}=(V,E)$ be any tree. It was shown by Valette (Theorem 1, \cite{V}) that, for any function $\psi : V \rightarrow \lbrack 0,1\rbrack$ satisfying the condition
$$  
\psi(x)\leq \frac{1}{\mathrm{deg}(x)},
$$
for all $x\in V$ (with the convention that $\psi(x)=0$ if $\mathrm{deg}(x)=\infty$), then, the kernel defined by
$$
\Psi(x,y)=\left\{
	\begin{array}{ll}
	0, & \textrm{if $x=y$}, \\
	d(x,y)-\frac{\psi(x)+\psi(y)}{2}, & \textrm{if $x\neq y$}, 
	\end{array}\right.
$$
is negative definite on $V$.\\
We address the following question. Let $G$ be a subgroup of $\mathrm{Aut}(\mathcal{T})$. When is $\Psi$ a $G$-invariant kernel? \\
We will answer this question in a special case.

\begin{Prop}\label{ValetteKernels}
Let $q\geq 2$ and let $G$ be a closed subgroup of $\mathrm{Aut}(\tk)$ acting transitively on both $V$ and $\bordtk$. Then, a kernel $\Psi$ defined as above is $G$-invariant if and only if the function $\psi$ used to construct $\Psi$ is constant.
\end{Prop}

\medskip

\noindent{\bf Proof : } Clearly, the kernel $\Psi$ is $G$-invariant if and only if the function $\psi$ satisfies the following condition
\begin{eqnarray}\label{invarianceCondition}
\psi(gx)+\psi(gy) &=& \psi(x)+\psi(y),
\end{eqnarray}
for any $x,y\in V$ and $g\in G$. Since the stabiliser of any vertex acts transitively on any sphere about any point, it is straightforward to see that $\psi$ has to take at most 2 values. Indeed, let us fix vertex $x_{0}$. Recall the standard bipartition of $V$ given by $V_{e}$ and $V_{o}$. The set $V_{e}$ (resp. $V_{o}$) consists of vertices which are at even (resp. odd) distance of $x_{0}$. Let $u,v$ be both in the same subset of the bipartition. Then, $d(u,v)$ is even and the median point of the geodesic $\lbrack u,v\rbrack$ is a certain vertex $z$. We can find $g\in \mathrm{stab}(z)$ sending $u$ on $v$. By condition (\ref{invarianceCondition}) we obtain
\begin{eqnarray*}
\psi(u)-\psi(v)
&=&
\psi(u)-\psi(gu) \\
&=&
\psi(gz)-\psi(z)\\
&=&
0,
\end{eqnarray*} 
which implies that $\psi(u)=\psi(v)$. To finish the proof, consider any geodesic segment of length $2$ formed by vertices $(v_{j})_{j=0}^{2}$. Since $G$ acts doubly transitively on $V$, we can find an element $g$ such that $gv_{j}=v_{j+1}$, for $j=0,1$. We observe that $d(v_{1},gv_{2})=2$ and this forces $\psi(gv_{2})=\psi(v_{1})$. Again, by condition (\ref{invarianceCondition}), we have
\begin{eqnarray*}
\psi(v_{0})
&=&
\frac{1}{2}\left(\psi(v_{0})+\psi(v_{2})\right) \\
&=&
\frac{1}{2}\left(\psi(gv_{0})+\psi(gv_{2})\right) \\
&=&
\frac{1}{2}\left(\psi(v_{1})+\psi(gv_{2})\right) \\
&=&
\psi(v_{1}),
\end{eqnarray*} 
and therefore, $\psi$ is constant. \hfill $\square$\\

\smallskip
 
\begin{Cor}
Let $G$ and $\Psi$ be as in Proposition \ref{ValetteKernels}. The only negative type functions on $G$ induced by the negative type kernels $\Psi$ are of the form
$$
g\mapsto d(x_{0},gx_{0})-\alpha,
$$
for some constant $\alpha\in\lbrack 0,\frac{1}{q+1}\rbrack$.
\end{Cor}

These functions are not pure in general, as we observe from the last three results. \\
Let $q\geq 2$ be a finite integer. As the decomposition of the Haagerup cocycle in Section \ref{decHaag} shows, the natural representation $\pi$ of a group $G \leq \textrm{Aut}(\tk)$ acting on $\ell_{\mathrm{alt}}^{2}(\mathbb{E})$ is reducible. Hence, we obtain the following corollary.

\begin{Cor}\label{distanceIsReducible}
The negative type function $g \mapsto d(x_{0},gx_{0})$ is not pure on $G$.
\end{Cor}

\medskip

We note that the result here above crucially depends on the fact that de degree of the tree is finite. Indeed, for the sake of completeness, we prove the following well-known Proposition.

\begin{Prop}
Let $q\geq 2$ be an integer (possibly infinite). Let $G$ be a closed non-compact subgroup of $\mathrm{Aut}(\tk)$. Suppose that $G$ acts transitively on the vertices and on the boundary $\bordtk$. Fix a vertex $x_{0}$, and set $|g|:=d(gx_{0},x_{0})$, for $g\in G$. Then, the function $g \mapsto |g|$ is pure in $\mathrm{CL}(G)$ if and only if $q=\infty$.
\end{Prop}
\noindent {\bf Proof : } Let us show that, in the case $q=\infty$, then the representation $\pi$ acting on $\ell_{alt}^{2}(\mathbb{E})$ is irreducible. Let us fix a geometric edge $a=\{a_{0},a_{1}\}\in E$. Firstly, we note that $\pi$ is equivalent to the quasi-regular representation $\lambda_{G/G_{a}}$, where $G_{a}$ is the stabiliser of $a$. By a theorem of Mackey (see Theorem 2.1 in \cite{BH}), we need to show that the commensurator of $G_{a}$ in $G$ is exactly $G_{a}$. Recall that the commensurator of $G_{a}$ in $G$, denoted by $\mathrm{Com}_{G}(G_{a})$, is the set of elements $g\in G$ such that the subgroup $G_{a}\cap G_{a'}$ has finite index in both $G_{a}$ and $G_{a'}$, where $a'$ is the edge satisfying $a'=ga$. Clearly, $G_{a}$ is contained in $\mathrm{Com}_{G}(G_{a})$.\\ 
To prove the other inclusion, let $g\in G\setminus G_{a}$ and set $a'=ga$ and $a_{i}'=ga_{i}$, for $i=0,1$. We will see that $G_{a}\cap G_{a'}$ has not finite index in $G_{a}$. We can suppose that the geodesic $\lbrack a_{0},a_{0}'\rbrack$ is contained in the geodesic $\lbrack a_{1},a_{1}'\rbrack$. Since the tree is of infinite degree, then, for all $k\geq 2$, there exists $a'_{k}\in V$ such that $a'_{0}\sim a'_{k}$ and $a'_{k}\neq a'_{j}$, for all $j\geq k-1$. Using the transitivity of the action on $\partial\mathcal{T}_{\infty}$, for every $k$, we can find $\tilde{g}_{k}\in G_{a_{0}}$ sending $a'$ to the edge $\{a'_{0},a'_{k}\}$. It is easy to see that $\tilde{g}_{k}\in G_{a}$, for all $k$, and that the cosets $\tilde{g}_{k}\left(G_{a}\cap G_{a'}\right)$ are pairwise different. This ends the proof. \hfill $\square$

\medskip 
 
We summarise the content of this section with a classification of pure elements of $\mathrm{CL}(G)$ when $G$ acts transitively on both $\tk$ and $\bordtk$.

\begin{Cor}
Let $G$ be a closed non-compact subgroup of $\mathrm{Aut}(\tk)$, with $q\geq 2$. Suppose that $G$ acts transitively on the vertices and on the boundary $\bordtk$. Let $\psi$ be a function conditionally of negative type on $G$. Suppose that $\psi$ is pure in $\mathrm{CL}(G)$ and that it vanishes on the stabiliser of some vertex $x_{0}$. We have the following alternative :
\begin{enumerate}
 	\item The function $\psi$ is bounded on $G$ and then it is of the form 
 	$$
 	\psi(g)=\|\xi\|^{2}-\langle \pi_{\psi}(g)\xi ,\xi\rangle ,
 	$$
 where $\pi_{\psi}$ is the irreducible unitary representation associated with $\psi$ via the GNS construction, and $\xi$ is a $\pi_{\psi}(G_{x_{0}})$-fixed vector (which is unique, up to scalar multiplication).
 	\item The function $\psi$ is unbounded and then it is of the form
 	$$
 	\psi(g)=C\left(|g|+\frac{2q}{q^{2}-1}(q^{-|g|}-1)\right),
 	$$
 	where $|g|:=d(gx_{0},x_{0})$ and $C$ is a positive constant.
\end{enumerate}  
\end{Cor}
\noindent {\bf Proof : } The only thing left to prove is the first claim. It is a general fact that $\psi(g)=\|(\pi_{\psi}(g)-1)\xi\|^{2}$, for some $G_{x_{0}}$-invariant vector $\xi\in\mathcal{H}_{\psi}$, using the GNS construction. By developing the norm and remarking that the coefficient $\langle \pi_{\psi}(\cdot)\xi,\xi\rangle$ is real-valued, we get $\psi(g)=2\|\xi\|^{2}-2\langle \pi_{\psi}(\cdot)\xi,\xi\rangle$. Finally, since the representation $\pi_{\psi}$ is spherical, then the space of $G_{x_{0}}$-invariant vectors has dimension one, and the result follows. \hfill $\square$
 
 \medskip


\bigskip

Authors addresses:
\medskip

\noindent
TU Dresden\\
Fachrichtung Mathematik\\
Institut für Geometrie\\
01062 Dresden

\begin{verbatim}antoine.gournay@tu-dresden.de
\end{verbatim}

\medskip

\noindent
Institut de Mathématiques - Unimail\\
11 Rue Emile Argand\\
CH-2000 Neuchâtel\\
Switzerland

\begin{verbatim}pierre-nicolas.jolissaint@unine.ch
\end{verbatim}

\end{document}